\documentclass[reqno, 10pt]{amsart}
\usepackage{graphicx}
\usepackage[dvips]{epsfig}
\usepackage{bm}
\usepackage{mathrsfs}

 \makeatletter
\@namedef{subjclassname@2020}{%
  \textup{2020} Mathematics Subject Classification}
\makeatother

\newtheorem{lemma}{Lemma}[section]
\newtheorem{proposition}[lemma]{Proposition}%[setcion]
\newtheorem{remark}[lemma]{Remark}%[section]
\newtheorem{problem}[lemma]{Problem}%[section]
\newtheorem{theorem}[lemma]{Theorem}%[section]

\newcommand \gam{\gamma}
\newcommand \R{\mathbb{R}}
\newcommand \Om{\Omega}
\newcommand \der{\partial}
\newcommand \vphi{\varphi}

\newcommand \mcl{\mathcal}
\newcommand \Gam{\Gamma}
\newcommand \alp{\alpha}
\newcommand \tx{\text}
\newcommand \til{\tilde}

\newcommand \ol{\overline}
\newcommand \eps{\varepsilon}

\newcommand \om{\omega}
\newcommand \rx{{\rm{\bf x}}}

\newcommand \Gamen{\Gamma_{0}}
\newcommand \Gamw{\Gamma_w}
\newcommand \Gamex{\Gam_L}

\newcommand \trho{\tilde{\rho}}
\newcommand \bvphi{\bar{\vphi}}
\newcommand \bPhi{\bar{\Phi}}
\newcommand \vrho{\varrho}
\newcommand \rhos{\vrho_s}
\newcommand \us{u_s}

\newcommand \tpsi{\til{\psi}}
\newcommand \tPsi{\til{\Psi}}
\numberwithin{figure}{section}

\numberwithin{equation}{section}

\allowdisplaybreaks
\begin{document}

\title[3-D Supersonic solutions to the Euler-Poisson system]{Three-dimensional Supersonic flows of Euler-Poisson system for potential flow}

%\author{Bae}
\author{Myoungjean Bae}
\address{Department of Mathematical Sciences, KAIST, 291 Daehak-Ro, Yuseong-Gu, Daejeon, South Korea 34141; Korea Institute for Advanced Study, 85 Hoegiro, Dongdaemun-Gu, Seoul 130-722, Republic of Korea}
\email{mjbae@kaist.ac.kr}
%\thanks{The first author was supported in part by NSF Grant \#000000.}

%\author{Park}
\author{Hyangdong Park}
\address{Center for Mathematical Analysis and Computation (CMAC), Yonsei University, 50 Yonsei-Ro, Seodaemun-Gu, Seoul 03722, Republic of Korea}
\email{hyangdong.park@yonsei.ac.kr}



 %35M32 (2010-now) Boundary value problems for mixed-type systems of PDEs
 %35L72 (2010-now) Second-order quasilinear hyperbolic equations
% 35J66 (2010-now) Nonlinear boundary value problems for nonlinear elliptic equations
 %35G60: Boundary value problems for systems of nonlinear higher-order PDEs
%35J47:  Second-order elliptic systems
%35J57: Boundary value problems for second-order elliptic systems
%35J66: Nonlinear boundary value problems for nonlinear elliptic equations
%35M10: Equations of mixed type
%35M30: Systems of mixed type????
%35Q31: Euler equations
%76N10: Compressible fluids and gas dynamics, general - Existence, uniqueness, and regularity theory
%76N10 (1980-now) Existence, uniqueness, and regularity theory for compressible fluids and gas dynamics [See also 35L60, 35L65, 35Q30]
%76J20 (1991-now) Supersonic flows

\date{\today}

\begin{abstract}
We prove the unique existence of supersonic solutions of the Euler-Poisson system for potential flow in a three-dimensional rectangular cylinder when prescribing the velocity and the strength of electric field at the entrance. Overall, the main framework is similar to \cite{BDXX}, but there are several technical differences to be taken care of vary carefully. And, it is our main goal to treat all the technical differences occurring when one considers a three dimensional supersonic solution of the steady Euler-Poisson system.
\end{abstract}

\maketitle

%\tableofcontents

%%%%%%%%%%%%%%%%%%%%%%%%%%%%%%%%%%%%%%
%%%%%%%%%%%%%%%%%%%%%%%%%%%%%%%%%%%%%%
\section{Introduction and the main problem}

Given a positive function $b({\bf x})$ representing the density of fixed, positively charged background ions, the steady Euler-Poisson system
\begin{equation}\label{E-S}
\left\{
\begin{split}
&\mbox{div}(\rho{\bf u})=0\\
&\mbox{div}(\rho{\bf u}\otimes{\bf u}+p\,{\mathbb I}_n)=\rho\nabla\Phi
%\quad(\mathbb{I}\mbox{ : $3\times 3$ identity matrix})
\\
&\mbox{div}\left(\rho(\mathcal{E}+\frac{p}{\rho}){\bf u}\right)=\rho{\bf u}\cdot\nabla\Phi \\
&\Delta\Phi=\rho-b
\end{split}
\right.
\end{equation}
yields a $n$-dimensional hydrodynamic model of semiconductor devices or plasmas, where ${\mathbb I}_n$ represents the $n\times n$ identity matrix.
%the motion of electrons governed by self-generated electric field in macroscopic scale.
Here, the functions $\rho=\rho({\bf x})$, $\displaystyle{{\bf u}=\sum_{j=1}^n u_j({\bf x})\hat{\bf e}_j}$, $p=p({\bf x})$, and $\mathcal{E}=\mathcal{E}({\bf x})$ represent the macroscopic particle electron density, velocity, pressure, and the total energy, respectively, at ${\bf x}\in\mathbb{R}^n$.
The function $\Phi=\Phi({\bf x})$ is the electric potential generated by the Coulomb force of particles.
%And, the function $b=b({\bf x})>0$ represents the density of fixed, positively charged background ions.

The ultimate goal of this work is to lay the foundation for the study of three dimensional supersonic flow of Euler-Poisson system in the regime of Helmholtz decomposition. And, we take the very first step in this paper by studying a simple problem of three dimensional supersonic solution to a potential flow model of steady Euler-Poisson system.

The potential flow model that we consider is derived from \eqref{E-S} under the assumption of ${\bf u}=\nabla\vphi$ for a scalar function $\vphi$, and under the constitutive relations of
\begin{equation*}
  p=\rho^{\gam},\quad {\mcl{E}}+\frac{p}{\rho}=\frac 12 |{\bf u}|^2+i(\rho)
\end{equation*}
for $i(\rho)$ given by
\begin{equation*}
  i(\rho)=\begin{cases}
  \frac{\gam(\rho^{\gam-1}-1)}{\gam-1}&\mbox{for $\gam>1$},\\
  \ln \rho &\mbox{for $\gam=1$}.
  \end{cases}
\end{equation*}
The case of $\gam>1$ concerns the ideal polytropic gas, and the case of $\gam=1$ concerns the isothermal gas. Under all the conditions stated in the above, if $\rho>0$ is assumed additionally, then the steady Euler-Poisson system \eqref{E-S} is simplified as follows:
\begin{equation}
\label{EP-pt-1}
  \begin{cases}
  &{\rm{div}} \left(\rho \nabla\vphi\right)=0\\
  &\nabla (\frac 12|\nabla\vphi|^2+i(\rho)-\Phi)=0\\
  &\Delta \Phi=\rho-b
  \end{cases}\quad.
\end{equation}
We replace the second equation in the above by
\begin{equation*}
 \frac 12|\nabla\vphi|^2+i(\rho)-\Phi=0,
\end{equation*}
and solve it for $\rho$ to obtain that
\begin{equation}
\label{definition-density}
  \rho=\trho(\Phi, \nabla\vphi):=\begin{cases}
  \left(1+\frac{\gam-1}{\gam}(\Phi-\frac 12|\nabla\vphi|^2)\right)^{\frac{1}{\gam-1}}&\mbox{for $\gam>1$}\\
  e^{(\Phi-\frac 12|\nabla\vphi|^2)}&\mbox{for $\gam=1$}
  \end{cases}\quad.
\end{equation}
So, for a fixed $\gam\ge 1$, we can further simplify the system \eqref{EP-pt-1} as
\begin{equation}
\label{EP-pt-final}
\begin{cases}
{\rm div}\left(\trho(\Phi, \nabla \vphi)\nabla\vphi\right)=0\\
\Delta \Phi=\trho(\Phi, \nabla\vphi)-b
\end{cases}\quad.
\end{equation}
As a second order quasi-linear system for $(\vphi, \Phi)$, the system \eqref{EP-pt-final} is \emph{a hyperbolic-elliptic mixed system} if and only if $|\nabla\vphi|^2> \gam \trho^{\gam-1}(\Phi, \nabla\vphi)$, and is \emph{an elliptic system} if and only if $|\nabla\vphi|^2< \gam \trho^{\gam-1}(\Phi, \nabla\vphi)$. As an analogy of the case of steady Euler system, the case of $|\nabla\vphi|^2> \gam \trho^{\gam-1}(\Phi, \nabla\vphi)$ corresponds to {\emph{a supersonic flow}}, and the case of $|\nabla\vphi|^2< \gam \trho^{\gam-1}(\Phi, \nabla\vphi)$ corresponds to {\emph{a subsonic flow}}. In this paper, we are aimed to construct a three dimensional supersonic solution to \eqref{EP-pt-final}. Even though the main framework would be similar to \cite{BDXX}, there are several technical differences to be taken care of, and it is our main goal to treat all the technical differences occurring when one considers a three dimensional supersonic solution of the steady Euler-Poisson system. This work may seem as a simple extension of \cite{BDXX} by adding a few technical changes. But we must emphasize that this work is a keystone to investigate a three dimensional supersonic solution to \eqref{E-S} with a nonzero vorticity in the regime of Helmholtz decomposition. Furthermore, the system \eqref{EP-pt-final} is a natural example of a quasi-linear hyperbolic-elliptic mixed system of second order arising from the study of physical phenomena, therefore the main result obtained from this paper can provide a footstone in the study of a nonlinear hyperbolic-elliptic mixed type system of second order partial differential equations in a multidimensional domain. Moreover, we remove the extra compatibility conditions given in \cite{BDXX} in this paper.

For a fixed constant $L>0$, we define a nozzle $\Omega_L$ of length $L$ with a rectangular cross section $\mathcal{D}:=(-1,1)\times(-1,1)\subset\mathbb{R}^2$ by
\begin{equation}\label{def-N-cyl}
\Omega_L:=\left\{(x_1,x_2,x_3)\in\mathbb{R}^3:\mbox{ }  0<x_1<L,\, (x_2,x_3)\in \mathcal{D}\right\}.
\end{equation}
The boundary of $\Om_L$ consists of the entrance $\Gamen$, the wall $\Gamw$, and the exit $\Gamex$ given respectively by
% The wall $\Gamma_w$, entrance $\Gamma_0$, and the exit $\Gamma_L$ of $\Omega_L$ are defined as follows:
\begin{equation*}
\Gamma_0:=\{0\}\times \ol{\mcl{D}},\quad
\Gamma_w:=(0,L)\times\partial \mathcal{D},\quad \Gamma_L:=\{L\}\times \ol{\mcl{D}}.
\end{equation*}

Suppose that $(\bvphi, \bPhi)(x_1)$ solves \eqref{EP-pt-final} with $b=b_0$ for some constant $b_0>0$, and let us set
\begin{equation}
\label{def-1d-potentials}
  u(x_1):=\bvphi'(x_1),\quad E(x_1):=\bPhi'(x_1),\quad \vrho(x_1):=\trho(\bPhi, \bvphi')(x_1).
\end{equation}
Then it can be directly checked that $(u, E, \vrho)$ solves the following ODE system:
\begin{equation}
\label{EP-1d}
 \begin{cases}
 u'=\frac{u^{\gam}E}{u^{\gam+1}-\gam J_0^{\gam-1}}\\
 E'=\frac{J_0}{u}-b_0\\
 \varrho=\frac{J_0}{u}
 \end{cases}
\end{equation}
for some constant $J_0>0$ provided that $u>0$. Note that the system \eqref{EP-pt-final} is derived from \eqref{E-S} under the assumption of $\rho>0$. Given a constant $J_0>0$, it can be directly checked that the state of  $(u, E, \vrho)$ is {\emph{supersonic}} if and only if $\varrho<\rhos$, and {\emph{subsonic}} if and only if $\varrho>\rhos$ for $\rhos$ given by
\begin{equation*}
  \rhos=\left(\frac{J_0^2}{\gamma }\right)^{\frac{1}{\gamma+1}}.
\end{equation*}

Let us fix two constants $\gam\ge 1$ and $J_0>0$.  And, assume that the constant $b_0$ in \eqref{EP-1d} satisfies
\begin{equation}
\label{assumption-b0}
  0<b_0<\rhos.
\end{equation}
\begin{remark}
The assumption of \eqref{assumption-b0} yields a description of a class of solutions to \eqref{EP-1d}. But the rest of argument for a multidimensional supersonic solution given in this paper is still valid even if this assumption is removed.
\end{remark}

For two constants $E_0\in \R$ and $u_0$ fixed with
\begin{equation*}
  u_0>\frac{J_0}{\rhos}=:\us,
\end{equation*}
let us consider the initial value problem
\begin{equation}
\label{ivp}
  \begin{cases}
 u'=\frac{u^{\gam}E}{u^{\gam+1}-\gam J_0^{\gam-1}}\\
 E'=\frac{J_0}{u}-b_0
 %\varrho=\frac{J_0}{u}
 \end{cases}\quad\mbox{with}\quad
 \begin{cases}
u(0)=u_0\\
E(0)=E_0
 \end{cases}
\end{equation}
on an interval $[0, L]$ for $L>0$ to be specified later.
\begin{lemma}[One dimensional supersonic solutions \cite{BDXX, LuoXin}]\label{one-sol}

For any small constant $\delta>0$, there exists a constant $L_1>0$ depending on $(\gam, J_0, b_0, u_0, E_0)$ so that the initial value problem \eqref{ivp} has a unique smooth solution $(u, E)$ on $[0, L_1]$ with satisfying that
\begin{equation*}
  \inf_{x_1\in [0, L_1]} u(x_1)\ge \us+\delta.
\end{equation*}
In other words, the solution $(u, E)$ provides a supersonic flow on the interval $[0, L_1]$.
\end{lemma}
A proof of Lemma \ref{one-sol} can be easily given from \cite[Lemma 1.1]{BDXX} or \cite{LuoXin} so we skip it in this paper. Instead, we add a brief comment on how $L_1$ varies depending on $(u_0, E_0)$. Any $C^1$ solution $(u, E)$ to \eqref{ivp} satisfies
\begin{equation*}
\frac{1}{2}E^2(x_1)-H(u(x_1))\equiv\frac 12 E_0^2-H(u_0)
\end{equation*}
for
$H(u)=\int_{\us}^{u} \frac{b_0}{t^{\gam+1}}(t^{\gam+1}-\us^{\gam+1})
\left(\frac{J_0}{b_0}-t\right)dt$. If $\frac{1}{2}E_0^2-H(u_0)<0$ holds, then $(u, E)$ is periodic and we have $u>\us$ for all $x_1>0$, and this implies that $L_1$ in Lemma \ref{one-sol} can be arbitrarily large. If $\frac{1}{2}E_0^2-H(u_0)\ge 0$ holds, on the other hand, the analysis on \cite[Lemma 1.1]{BDXX} or \cite{LuoXin} shows that $L_1$ in Lemma \ref{one-sol} is finite. We point out that this property of $L_1$ is given due to the assumption of \eqref{assumption-b0}.

For the smooth supersonic solution $(u, E)$ to the initial value problem \eqref{ivp}, the vector field $(\bar{\vphi},\bar{\Phi})$ in \eqref{def-1d-potentials} given by
\begin{equation*}%\label{def-phi0}
\begin{split}
\bar{\vphi}(x_1):=\int_0^{x_1} {u}(t)\,dt,\quad
\bar{\Phi}(x_1):=\frac 12 u_0^2+i\left(\frac{J_0}{u_0}\right)+\int_0^{x_1} E(t)\,dt
\end{split}
\end{equation*}
provides a one dimensional solution to \eqref{EP-pt-final} in $\Om_L$ for any $L\in(0, L_1]$. Let us define
\begin{equation}\label{def-phi0}
\begin{split}
{\vphi}_0({\rx}):=\bar{\vphi}(x_1),\quad
{\Phi}_0(\rx):=\bar{\Phi}(x_1)\quad\tx{for ${\rx}=(x_1, x_2, x_3)\in \Om_L$}.
\end{split}
\end{equation}

%Now we state the problem and main results.
The main goal of this work is to construct a supersonic solution to \eqref{EP-pt-final} as a small perturbation of $(\vphi_0, \Phi_0)$ in $\Om_L$ by solving the problem stated in the following:
\begin{problem}\label{Problem}
Given functions $b\in C^2(\overline{\Omega_L})$, $u_{\rm en}\in C^3(\overline{\Gamma_0})$, $E_{\rm en}\in C^4(\overline{\Gamma_0})$, and $\Phi_{\rm ex}\in C^4(\overline{\Gamma_L})$, we prescribe the boundary conditions for $(\vphi, \Phi)$ as follows:
\begin{equation}\label{irro-bd}
	\left.\begin{split}
	\varphi=0,\quad \partial_1\varphi=u_{\rm en},\quad \partial_1\Phi=E_{\rm en}\quad&\mbox{on}\quad\Gamma_0,\\
	\der_{\bf n}\varphi=0,\quad \der_{\bf n}\Phi=0\quad&\mbox{on}\quad\Gamma_w,\\
	\Phi=\Phi_{\rm ex}\quad&\mbox{on}\quad \Gamma_L,
	\end{split}\right.
	\end{equation}
where ${\bf n}$ represents the inward unit normal on $\Gamma_w$.
And, let us set
\begin{equation*}
\begin{split}
\sigma(b,u_{\rm en},E_{\rm en},\Phi_{\rm ex}):=&\|b-b_0\|_{C^2(\overline{\Omega_L})}+\|u_{\rm en}-u_0\|_{C^3(\overline{\Gamma_0})}\\
&+\|E_{\rm en}-E_0\|_{C^4(\overline{\Gamma_0})}+\|\Phi_{\rm ex}-\Phi_0\|_{C^4(\overline{\Gamma_L})}.\end{split}
\end{equation*}

For $\sigma(b,u_{\rm en},E_{\rm en},\Phi_{\rm ex})$ sufficiently small, solve \eqref{EP-pt-final} in $\Om_L$ with the boundary conditions stated in \eqref{irro-bd}.

\end{problem}

\begin{remark}
\label{remark-corner-slip-bc}
For ${\rm x}\in \Gam_w\cap\{|x_2|=1\,\,\tx{and}\,\,|x_3|=1\}$, we interpret the slip boundary condition $\der_{\bf n} \vphi=0$ in the following sense:
\begin{equation*}
  \lim_{{\rm y\to \rm x}\atop {{\rm y}\in \Gam_w\setminus\{|x_2|=1\,\,\tx{and}\,\,|x_3|=1\}}} \der_{\bf n} \vphi({\rm y})=0.
\end{equation*}
And, the slip boundary condition $\der_{\bf n} \Phi=0$ on $\Gam_w$ is interpreted in the same sense. In other words, we expect for a solution $(\vphi,\Phi)$ of Problem \ref{Problem} to satisfy the conditions
\begin{equation*}
  \der_{x_k}\vphi=0,\quad\tx{and}\quad
  \der_{x_k}\Phi=0\quad\tx{on $\Gam_w\cap\{|x_2|=1\,\,\tx{and}\,\,|x_3|=1\}$ for $k=2,3$}.
\end{equation*}
\end{remark}

%%%%%%%%%%%%%%%%%%%%%%%%%%%%%%%%%%%%%

%%%%%%%%%%%%%%%%%%%%%%%%%%%%%%%%%%%%%%%%
\begin{theorem}\label{Irr-Thm}
Given a small constant $\delta>0$, let the constant $L_1$ be from Lemma \ref{one-sol}. Then, there exists a constant $L_1^*\in(0, L_1]$ depending on $\gam$, $J_0$, $u_0$, $E_0$, $b_0$, and $\delta$ so that the following properties hold:
For a fixed $L<L_1^*$, one can fix a small constant $\bar{\sigma}>0$ depending on $\gam$, $J_0$,  $u_0$, $E_0$, $b_0$, $\delta$, and $L$ so that if
\begin{equation*}
  \sigma(b,u_{\rm en},E_{\rm en},\Phi_{\rm ex})\le \bar{\sigma},
\end{equation*}
and if the following compatibility conditions holds:
\begin{equation}\label{223}
\begin{split}
\der_{{\bf n}}b=0&\quad\mbox{on $\Gam_w$},\\
\der_{\bf n}u_{\rm en}=0,\,\, \der^k_{{\bf n}}E_{\rm en}=0&\quad\mbox{on $\ol{\Gam_w}\cap \ol{\Gam_0}$ for $k=1,3$},\\
\der^k_{{\bf n}}\Phi_{\rm ex}=0
&\quad\mbox{on $\ol{\Gam_w}\cap \ol{\Gam_L}$ for $k=1,3$},
%\frac{db}{dx_2}=\frac{d^ku_{\rm en}}{dx_2^k}=\frac{d^kE_{\rm en}}{dx_2^k}=\frac{d^k\Phi_{\rm ex}}{dx_2^k}=0&\mbox{ on }\Gamma_w\cap\{x_2=\pm1\}\mbox{ for }\,k=1,3,\\
%\frac{db}{dx_3}=\frac{d^ku_{\rm en}}{dx_3^k}=\frac{d^kE_{\rm en}}{dx_3^k}=\frac{d^k\Phi_{\rm ex}}{dx_3^k}=0&\mbox{ on }\Gamma_w\cap\{x_3=\pm1\}\mbox{ for }\,k=1,3,
\end{split}
\end{equation}
then Problem \ref{Problem}  has a unique solution $(\varphi,\Phi)\in [H^4(\Omega_L)]^2$ that satisfies the estimate
\begin{equation}\label{Thm-est}
\|\varphi-\varphi_0\|_{H^4(\Omega_L)}+\|\Phi-\Phi_0\|_{H^4(\Omega_L)}\le C\sigma(b,u_{\rm en},E_{\rm en},\Phi_{\rm ex})\end{equation}
for a constant $C>0$ depending only on $(\gamma, J_0,b_0,u_0,E_0,\delta,L)$.

%Furthermore, it holds that $\Phi\in C^2(\Omega_L)$.
\end{theorem}

Before we proceed to prove Theorem \ref{Irr-Thm}, we first add a few remarks on Theorem \ref{Irr-Thm}.

\begin{remark}
\label{remark1-main theorem}
It follows from \eqref{Thm-est} and the generalized Sobolev inequality  that
%\begin{equation}\label{C-est}
%\|\varphi-\varphi_0\|_{C^{2,\frac{1}{2}}(\overline{\Omega_{L}})}+\|\Phi-\Phi_0\|_{C^{2,\frac{1}{2}}(\overline{\Omega_{L}})}\le C_{\sharp}\sigma(b,u_{\rm en},E_{\rm en},\Phi_{\rm ex})
%\end{equation}
%for a constant $C_{\sharp}>0$ depending only on $(\gamma, J_0,b_0,u_0,E_0,\delta,L)$. Therefore,
the solution $(\vphi, \Phi)$ in Theorem \ref{Irr-Thm} is a $C^2$ classical solution of Problem \ref{Problem}.

\end{remark}

\begin{remark}[The restriction on $L$ and the dependence of $\bar{\sigma}$ on $L$]
\label{remark2-main theorem}
For readers familiar with the study on three dimensional supersonic solution of steady Euler system, the additional restriction of $L\in (0, L_1^*]$ in Theorem \ref{Irr-Thm}
%, and the dependence of the smallness of $\sigma(b,u_{\rm en},E_{\rm en},\Phi_{\rm ex})$ on $L$
may seem strange or redundant even. We must emphasize that the restriction on $L$ is given in order to establish a priori $H^1$ energy estimate for hyperbolic-elliptic coupled system of second order. And, this is the main difference of the Euler-Poisson system from the Euler system. Namely, the equation $\Delta \Phi=\rho-b$ for the electric potential $\Phi$, which is elliptic, is coupled with the rest of the Euler-Poisson system in a nonlinear way, and it causes many technical challenges ranging from an appropriate choice of a solution space to developing a method to establish a priori energy estimate for solutions.
\end{remark}

\begin{remark}
\label{remark3-main theorem}
It is our ultimate goal
%The ultimate goal of this initiative work for multidimensional supersonic solutions to the Euler-Poisson system is
to achieve the well-posedness
%of Problem \ref{Problem} and
of a boundary value problem for a supersonic flow of full Euler-Poisson system with nonzero vorticity in a three or higher dimensional cylindrical domain with an arbitrary cross section. In a three or higher dimensional cylindrical domain with an arbitrary cross section, it is well known that one can construct a smooth supersonic solution to the steady Euler system for potential flow as long as appropriate compatibility conditions for boundary data are prescribed at Lipschitz corners of the cylindrical domain. But it is not the case for the steady Euler-Poisson system \eqref{E-S}(or \eqref{EP-pt-final}) mainly because the elliptic equation $\Delta \Phi=\rho-b$ for the electric potential $\Phi$ has a nonlocal effect to the rest of equations in \eqref{E-S} as the density $\rho$ is determined depending on $\Phi$. So it is difficult to find a supersonic solution to \eqref{E-S} or \eqref{EP-pt-final} with satisfying the compatibility conditions.
As an initiative work, we consider a potential flow in a rectangular cylinder in this paper. Supersonic flows of Euler-Poisson system with nonzero vorticity in a cylindrical domain with a more general cross sections will be studied in the forthcoming works in the regime of Helmholtz decomposition, for which the study of a potential flow model is crucial.

\end{remark}

We explain how to prove Theorem \ref{Irr-Thm} in the next section. Since the proof is more or less similar to the proof of \cite[Theorem 1.7]{BDXX} except that a priori $H^s$ estimate of approximate solutions for $s=1,\cdots, 4$ is more complicated, we provide an outline of a proof in Section \ref{subsection-outline} first, then we discuss further details on the technical differences in a priori $H^s$ estimates for the three dimensional case.
% from the two dimensional case studied in \cite{BDXX}.

%%%%%%%%%%%%%%%%%%%%%%%%%%%%%%%%%%%%%%%%%
%%%%%%%%%%%%%%%

\section{Proof of Theorem \ref{Irr-Thm}}
\subsection{Outline} \label{subsection-outline}
For $\til{\rho}$ given by \eqref{definition-density}, let us define $\til{c}$ by
\begin{equation}
\label{Sound}
  \til{c}(\Phi, \nabla\vphi):=\sqrt{\gam\til{\rho}^{\gam-1}(\Phi, \nabla\vphi)}.
\end{equation}
For $z\in\mathbb{R}$, ${\bf p}=(p_1,p_2,p_3)\in\mathbb{R}^3$, ${\bf q}=(q_1,q_2,q_3)\in\mathbb{R}^3$ with $z-\frac{1}{2}|{\bf q}|^2>0$, define
\begin{equation}\label{Aij}
\begin{split}
&A_{ij}(z,{\bf q}):=\frac{\til c^2(z,{\bf q})\delta_{ij}-q_iq_j}{\til c^2(z,{\bf q})-q_1^2}\quad\mbox{for }i,j=1,2,3,\\
&B(z,{\bf p},{\bf q}):=\frac{{\bf p}\cdot{\bf q}}{\til c^2(z,{\bf q})-q_1^2}.
\end{split}
\end{equation}
Here, $\delta_{ij}$ represents the Kronecker delta, that is, $\delta_{ij}=1$ for $i=j$, and $0$ otherwise.

As we seek a classical solution $(\vphi, \Phi)$ to \eqref{EP-pt-final} with satisfying \eqref{Thm-est} for $\sigma(b, u_{\rm en}, E_{\rm en}, \Phi_{\rm ex})$ sufficiently small, we expect that
\begin{equation*}
\til{c}^2(\Phi,\nabla\varphi)-(\partial_1\varphi)^2<0\quad\mbox{in}\quad \overline{\Omega_L}.
\end{equation*}
So we rewrite \eqref{EP-pt-final} as
\begin{equation}\label{First}
\left\{\begin{split}
	&\begin{split}
	\sum _{i,j=1}^3 A_{ij}(\Phi,\nabla\varphi)\partial_{ij}\varphi+B(\Phi,\nabla\Phi,\nabla\varphi)&=0,
	\end{split}\\
	&\Delta\Phi=\til{\rho}(\Phi,\nabla\varphi)-b.
\end{split}\right.
\end{equation}

Next, we further rewrite \eqref{First} in terms of the perturbations
$$(\psi, \Psi):=(\varphi-\varphi_0, \Phi-\Phi_0).$$
A lengthy but straightforward computation yields that
\begin{equation}
\label{Eqn-perturbation}
\left\{\begin{split}
&\sum_{i,j=1}^3 a_{ij}(\Psi,\nabla\psi)\partial_{ij}\psi
+\bar{a}_1\partial_1\psi+\bar{b}_1\partial_1\Psi+\bar{b}_2 \Psi=\mathfrak{f}_1(\Psi,\nabla\Psi,\nabla\psi),\\
&\Delta\Psi-\bar{h}_1\Psi-\bar{h}_2\der_1\psi=\mathfrak{f}_2(\Psi,\nabla\psi),
\end{split}\right.
\end{equation}
where $a_{ij}$ for $1\le i,j\le 3$, $\bar{a}_1$, $\bar{b}_k$, $\bar{h}_k$ and $\mathfrak{f}_k$ for $k=1,2$ are given as follows:
\begin{align}
\label{def-aij}
&\begin{cases}
a_{ij}(\xi,{\bm \zeta}):=A_{ij}(\xi+\Phi_0,{\bm \zeta}+\nabla\varphi_0),\\
(\bar{a}_1,\bar{b}_1,\bar{b}_2):=
\left.(\partial_{q_1},\partial_{p_1},\partial_z)B(z,{\bf p},{\bf q})\right|_{(z,{\bf p},{\bf q})=(\Phi_0,\nabla\Phi_0,\nabla\varphi_0)},\\
(\bar{h}_1,\bar{h}_2):=\left.(\partial_z,\partial_{q_1})\til{\rho}(z,{\bf q})\right|_{(z,{\bf q})=(\Phi_0,\nabla\varphi_0)},
\end{cases}\\
\label{def-f1}
&\mathfrak{f}_1(\xi,{\bm \eta},{\bm \zeta}):=-B(\Phi_0+t\xi,\nabla\Phi_0+t{\bm \eta},\nabla\varphi_0+t{\bm \zeta})|_{t=0}^1
+\bar{a}_1\zeta_1+\bar{b}_1\eta_1+\bar{b}_2\xi,\\
\label{def-f2}
&\mathfrak{f}_2(\xi,{\bm \zeta}):=
\til{\rho}(\Phi_0+t\xi,\nabla\varphi_0+t{\bm \zeta})|_{t=0}^1
-(b-b_0)-\bar{h}_1\xi-\bar{h}_2\zeta_1
\end{align}
for $z, \xi\in \R$, ${\bf p}, {\bf q}, {\bm \eta}, {\bm \zeta}\in \R^3$.

In the above, $z$ and $\xi$ are symbols for functions given in $\ol{\Om_L}$. And, ${\bf p}$, ${\bf q}$, ${\bm \eta}$ and ${\bm \zeta}$ are symbols for gradients of the functions. So, at each ${\rm x}\in \ol{\Om_L}$, we presume that $(\Phi_0, \nabla\Phi_0, \nabla\vphi_0)$ are evaluated at ${\rm x}$ thus  $a_{ij}$, $\bar{a}_1$, $\bar{b}_k$, $\bar{h}_k$ and $\mathfrak{f}_k$ are functions given in $\Om_L$ once $(\xi, {\bm \eta}, {\bm \zeta})$ are fixed as $(\Psi, \nabla\Psi, \nabla\psi)$.

Fix two functions $\til{\psi}$ and $\til{\Psi}$ in $\ol{\Om_L}$.
For later use, we define two linear operators $\mathcal{L}_1$ and $\mathcal{L}_2$ associated with $(\tpsi, \tPsi)$ by
\begin{equation}\label{def-L2}
\begin{split}
&\mathcal{L}_1^{(\tpsi, \tPsi)}(\psi,\Psi)
:=\sum _{i,j=1}^3a_{ij}(\tPsi,\nabla\tpsi)\partial_{ij}\psi
+\bar{a}_1\partial_1\psi+\bar{b}_1\partial_1\Psi+\bar{b}_2 \Psi,\\
&\mathcal{L}_2(\psi,\Psi)
:=\Delta\Psi-\bar{h}_1\Psi-\bar{h}_2\partial_1\psi.
\end{split}
\end{equation}
We will further specify a class of $(\tpsi, \tPsi)$ so that $\mcl{L}_1^{(\tpsi, \tPsi)}$ is well defined in $\Om_L$.

Theorem \ref{Irr-Thm} can be proved by solving the following boundary value problem:
\begin{equation}\label{nlbvp-perturbation}
\left\{\begin{split}
\mathcal{L}_1^{(\psi, \Psi)}(\psi,\Psi)=\mathfrak{f}_1(\Psi,\nabla\Psi,\nabla\psi),\quad
\mathcal{L}_2(\psi,\Psi)=\mathfrak{f}_2(\Psi,\nabla\psi)\quad&\mbox{in}\quad\Omega_L,\\
(\psi,\partial_1\psi)=(0,u_{\rm en}-u_0)=:(0,g_1),\quad\partial_1\Psi=E_{\rm en}-E_0=:g_2\quad&\mbox{on}\quad\Gamma_0,\\
\Psi=\Phi_{\rm ex}-\Phi_0(L,\cdot)=:\Psi_{\rm ex}\quad&\mbox{on}\quad\Gamma_L,\\
\der_{\bf n}\psi=0,\quad \der_{\bf n}\Psi=0\quad &\mbox{on}\quad\Gamma_w.
\end{split}\right.
\end{equation}

Without any technical details, we first outline how to solve \eqref{nlbvp-perturbation}.
\medskip

Fix a small constant $\delta>0$, and let the constant $L_1$ be from Lemma \ref{one-sol}. For a small constant $\eps>0$, and a constant $L\in(0, L_1]$ to be specified later, let us define a set $\mcl{I}_{\eps, L}$ by
\begin{equation*}
\mathcal{I}_{\eps,L}:=\left\{\phi \in H^4(\Omega_L):\,
\|\phi\|_{H^4(\Omega_L)}\le\eps,\,\der_{\bf n}\phi=0\mbox{ on }\Gamma_w\right\},
\end{equation*}
then we define an iteration set
\begin{equation}\label{Ite-set}
\mathcal{J}_{\eps, L}:=\mathcal{I}_{\eps,L}\times\mathcal{I}_{\eps,L}.
\end{equation}

Fix $(\phi, W)\in \mcl{J}_{\eps, L}$, and we consider the following linear boundary value problem associated with $(\phi, W)$:
\begin{equation}\label{lin-pro2}
\left\{\begin{split}
\mathcal{L}_1^{(\phi, W)}(\psi,\Psi)
=\mathfrak{f}_1(W,\nabla W,\nabla \phi),\quad \mathcal{L}_2(\psi,\Psi)=\mathfrak{f}_2(W,\nabla \phi)\quad&\mbox{in}\quad\Omega_L,\\
(\psi,\partial_1\psi)=(0,g_1),\quad\partial_1\Psi=g_2\quad&\mbox{on}\quad\Gamma_0,\\
\Psi=\Psi_{\rm ex}\quad&\mbox{on}\quad\Gamma_L,\\
\der_{\bf n}\psi=0,\quad \der_{\bf n}\Psi=0\quad &\mbox{on}\quad\Gamma_w.
\end{split}\right.
\end{equation}

\begin{proposition}[Analogy of {\cite[Proposition 2.6]{BDXX}}]\label{proposition-main} Given a small constant $\delta>0$, there exists a constant $L_1^*\in(0, L_1]$ depending on $\gam$, $J_0$, $u_0$, $E_0$, $b_0$, and $\delta$ so that the following properties hold: For a fixed $L<L_1^*$, one can fix a small constant $\bar{\eps}>0$ depending on $\gam$, $J_0$,  $u_0$, $E_0$, $b_0$, $\delta$, and L so that if $\mcl{J}_{\eps, L}$ is given by \eqref{Ite-set} for $\eps$ satisfying $0<\eps\le \bar{\eps}$, then, for each $(\phi, W)\in \mcl{J}_{\eps, L}$, and any given functions $f_1\in H^3(\Om_L)$, $f_2\in H^3(\Om_L)$ and $g\in C^3(\ol{\Gam_0})$ with the compatibility conditions
\begin{equation}
\label{compatibility-original}
  \der_{\bf n} f_1=0,\quad \der_{\bf n} f_2=0\quad\mbox{on $\Gam_w$},\quad\tx{and}\quad \der_{\bf n}g=0\quad\tx{on $\ol{\Gam_0}\cap \ol{\Gam_w}$},
\end{equation}
the following linear boundary value problem
\begin{equation}\label{lin-bd}
\left\{\begin{split}
\mathcal{L}_1^{(\phi, W)}(\psi,\Psi)
=f_1,\quad
\mathcal{L}_2(\psi,{\Psi})=f_2\quad&\mbox{in}\quad\Omega_L,\\
(\psi,\partial_1\psi)=(0,g),\quad\partial_1{\Psi}=0\quad&\mbox{on}\quad\Gamma_0,\\
{\Psi}=0\quad&\mbox{on}\quad\Gamma_L,\\
\der_{\bf n}\psi=0,\quad \der_{\bf n}{\Psi}=0\quad &\mbox{on}\quad\Gamma_w
\end{split}\right.
\end{equation}
has a unique solution $(\psi,{\Psi})\in[H^4(\Omega_L)]^2$ that satisfies the estimate
\begin{equation}\label{hat-est}
\begin{split}
&\|\psi\|_{H^4(\Omega_L)}\le C\left(\|f_1\|_{H^3(\Omega_L)}
+\|f_2\|_{H^2(\Omega_L)}+\|g\|_{C^3(\overline{\Gamma_0})}\right),\\
&\|{\Psi}\|_{H^4(\Omega_L)}\le C\left(\|f_1\|_{H^2(\Omega_L)}+\|f_2\|_{H^2(\Omega_L)}+\|g\|_{C^2(\overline{\Gamma_0})}\right),
\end{split}
\end{equation}
for the estimate constant $C>0$ depending only on $(\gamma, J_0,  b_0, u_0, E_0, \delta, L)$.
\end{proposition}
Note that it is actually sufficient to assume that $f_2\in H^2(\Om_L)$ in Proposition \ref{proposition-main}. If it is assumed that $f_2$ is $H^3$ in $\Om_L$, then the compatibility condition $\der_{\bf  n}f_2=0$ on $\Gam_w$ should be considered in the sense of trace. Since $\mathfrak{f}_2(W, \nabla\phi)$ from \eqref{lin-pro2} is $H^3(\Om_L)$ for $(\phi, W)\in \mcl{J}_{\eps, L}$, we assume in Proposition \ref{proposition-main} that $f_2$ is $H^3$ in $\Om_L$ for simplicity.

The following lemma yields useful properties of the coefficients of the linear operator $\mcl{L}_1^{(\phi, W)}$ for $(\phi, W)\in \mcl{J}_{\eps, L}$.
\begin{lemma}\label{lem-nei}
Given a small constant $\delta>0$, let $L_1$ be from Lemma \ref{one-sol}.
\begin{itemize}
\item[(a)] For $a_{ij}(\xi, {\bm \zeta})$ given by \eqref{def-aij}, let us set
   $
\bar{a}_{ij}({\bf x}):=a_{ij}(0, {\bf 0})$  for $1\le i,j\le 3.$ Then, $\bar{a}_{ij}$ for $1 \le i,j\le 3$, $\bar{a}_1$, $\bar{b}_k$ and $\bar h_k$ for $k=1,2$ are smooth in $\overline{\Omega_L}$, and for each $l\in \mathbb{N}$, there exists a constant $M_l>0$ depending only on $(\gamma, J_0,  b_0, u_0, E_0, \delta, l)$  to satisfy
\begin{equation*}
\|(\bar{a}_{ij},\bar{a}_1,\bar{b}_k,\bar{h}_k)\|_{C^l(\overline{\Omega_L})}\le M_l.
\end{equation*}
Furthermore, there exists a constant $\mu_1\in(0,1)$ depending only on $\gamma$, $J_0$,  $b_0$, $u_0$, $E_0$ and $\delta$ so that the inequalities
\begin{equation*}
\mu_1\le \bar{a}_{11}, -\bar{a}_{jj}({\bf x})\le \frac{1}{\mu_1}%\quad\mbox{for}\quad {\bf x}\in.
\end{equation*}
hold in $\overline{\Omega_L}$ for $j=2,3$. In fact, $\bar a_{11}\equiv 1$ in $\ol{\Om_L}$.

\item[(b)] There exists a constant $\eps_1>0$ so that if $\eps\le \eps_1$ in \eqref{Ite-set}, then for each $(\phi, W)\in \mcl{J}_{\eps, L}$, the coefficients $a_{ij}(W, \nabla\phi)$ given by \eqref{def-aij} for $1\le i,j \le 3$ satisfy
    \begin{equation*}
      a_{ij}(W, \nabla\phi)=a_{ji}(W, \nabla\phi)\quad\tx{for $1\le i,j\le 3$},
    \end{equation*}
    \begin{equation*}
      \sup_{1\le i,j\le 3}\|a_{ij}(W, \nabla\phi)-\bar{a}_{ij}\|_{H^3(\Om_L)}
      \le C\eps,
    \end{equation*}
    for some constant $C>0$.
    Moreover, the matrix $-[a_{ij}(W, \nabla\phi)]_{i,j=2}^3$ is positive definite in $\Om_L$ with
    \begin{equation*}
      \frac{\mu_1}{2}\le -[a_{ij}(W, \nabla\phi)]_{i,j=2}^3\le \frac{2}{\mu_1}\quad\tx{in $\Om_L$}.
    \end{equation*}
    Here, the constants $\eps_1$ and $C$ can be fixed depending only on $\gamma$, $J_0$, $b_0$, $u_0$, $E_0$ and $\delta$.

\end{itemize}
\end{lemma}

Once Proposition \ref{proposition-main} is proved, Theorem \ref{Irr-Thm} can be proved by employing a method of iteration. One can refer to \cite[Proof of Theorem 1.6]{BDXX} for further details on how to apply Proposition \ref{proposition-main} in implementing the iteration process. In this paper, we focus on how to prove Proposition \ref{proposition-main}. In particular, we intend to point out technical difficulties that one may confront in dealing with a three dimensional cylindrical domain with an arbitrary cross section, and how to resolve the issues for a special case at the moment. More details are discussed in the next section.

\subsection{A remark on a Proposition \ref{proposition-main}}
\label{subsection-further-discussion}
Fix a small constant $\delta>0$, and let $L_1>0$ be from Lemma \ref{one-sol}. And, let us assume that $0<L\le L_1$.
Given $(\phi, W)\in \mcl{J}_{\eps, L}$, let $a_{ij}(W, \nabla\phi)$ be given by \eqref{def-aij}. To simplify notations, we let $a_{ij}$ denote $a_{ij}(W, \nabla\phi)$ unless otherwise specified. And, we let $\mcl{L}_1$ denote the linear operator $\mcl{L}_1^{(\phi, W)}$ given by \eqref{def-L2}. Take functions $f_1\in H^3(\Om_L)$, $f_2\in H^3(\Om_L)$ and $g\in C^3(\ol{\Gam_0})$ with the compatibility conditions
\begin{equation}
\label{compatibility}
  \der_{\bf n} f_1=0,\quad \der_{\bf n} f_2=0\quad\mbox{on $\Gam_w$},\quad\tx{and}\quad \der_{\bf n}g=0\quad\tx{on $\ol{\Gam_0}\cap \ol{\Gam_w}$}.
\end{equation}
And, we consider the following linear boundary value problem for $(\psi, \Psi)$:
\begin{equation}\label{lin-bd-simp}
\left\{\begin{split}
\mathcal{L}_1(\psi,\Psi)
=f_1,\quad
\mathcal{L}_2(\psi,{\Psi})=f_2\quad&\mbox{in}\quad\Omega_L,\\
(\psi,\partial_1\psi)=(0,g),\quad\partial_1{\Psi}=0\quad&\mbox{on}\quad\Gamma_0,\\
{\Psi}=0\quad&\mbox{on}\quad\Gamma_L,\\
\der_{\bf n}\psi=0,\quad\der_{\bf n}{\Psi}=0\quad &\mbox{on}\quad\Gamma_w.
\end{split}\right.
\end{equation}
The linear operators $\mcl{L}_1$ and $\mcl{L}_2$  are second order differential operators with respect to $\psi$ and $\Psi$, respectively, therefore we determine their types depending on their second order derivative terms. According to Lemma \ref{lem-nei}, the operator $\mcl{L}_1$ is hyperbolic type with respect to $\psi$ if $\eps$ is fixed sufficiently small in the definition of the iteration set $\mcl{J}_{\eps, L}$ so the following lemma is obtained:
\begin{lemma}\label{lemma-L1-hyperbolic-type}
There exists a constant $\eps_2\in(0,\eps_1]$ for $\eps_1$ from the statement (b) of Lemma \ref{lem-nei} so that if $\eps\le \eps_2$ in the definition \eqref{Ite-set} of the iteration set $\mcl{J}_{\eps, L}$, then the linear boundary value problem \eqref{lin-bd-simp} is a hyperbolic-elliptic coupled system.
\end{lemma}

In order to solve the linear boundary value problem \eqref{lin-bd-simp}, we establish a priori $H^4$ estimates of $(\psi, \Psi)$ and apply the method of Galerkin's approximations. Note that we seek for a $H^4$ solution of \eqref{lin-bd-simp} as we intend to find a $C^2$ classical solution of \eqref{lin-bd-simp}. Also, we discover that seeking for a $H^3$ solution of \eqref{lin-bd-simp} is not possible at least by applying the method developed in \cite{BDXX} and this paper because we heavily use the generalized H\"{o}lder inequality and Sobolev inequality in a delicate way, and this approach fails if we seek for a $H^3$ solution. %This is a surprising discovery that we find in this work.

The main difference of this work from \cite{BDXX} is that we find a general structure of the operator $\mcl{L}_1$ so that a priori $H^4$ estimates of solutions to \eqref{lin-bd-simp} can be established in any dimensional domain. What we mean by the general structure is not only the hyperbolicity of $\mcl{L}_1$ but also a particular structure of the coefficients $a_{ij}$. The quasi-nonlinearity of the first equation in \eqref{First} with respect to $\vphi$ leads to the fact that $a_{ij}$ depends on $\nabla\phi$ for $(\phi,W)\in \mcl{J}_{\eps, L}$ satisfying the slip boundary condition $\der_{\bf n}\phi=0$ on $\Gam_w$. Thanks to this useful structure, the computation for a priori $H^1$ energy estimate given in \cite{BDXX} can be extended to any dimensional domain by a careful treatment.

\subsection{A brief description on how to prove Proposition \ref{proposition-main}}\label{subsection-more-details}
The rest of the paper is devoted to provide more technical details which must be clarified for the three dimensional domain $\Om_L$ of our consideration so that they can be easily cooperated with the results from \cite{BDXX} in order to prove Proposition \ref{proposition-main}. Once the details are given in this paper, then readers can refer to \cite{BDXX} to complete to prove Proposition \ref{proposition-main} thus prove Theorem \ref{Irr-Thm}, the main result of this paper.

\subsubsection{Approximated linear boundary value problem}\label{subsubsection-approximation} To avoid any technical constraint of weak regularity in implementing various computations towards $H^4$ a priori estimates, we first approximate the coefficients $a_{ij}$ by smooth functions via the method of extensions. In order to keep a useful property of $a_{ij}$, we provide smooth extensions of $(\phi, W)\in \mcl{J}_{\eps, L}$ first, then use them to have smooth approximations of $a_{ij}$. Even though the following procedure is somewhat typical, we put details as they contribute to simplify the definition \eqref{Ite-set} of the iteration set $\mcl{J}_{\eps, L}$ more than the one given in \cite[(2.19)]{BDXX}. Moreover, the extension method introduced in the below can be applied to a multidimensional cylindrical domain with a general cross section with a Lipschitz boundary.

Given $\phi\in \mcl{I}_{\eps, L}$ we define its first extension $\phi_{(1)}$ onto $\Om_L^{(1)}:=\{{\rm x}\in \R^3:0<x_1<L,\, |x_2|<\frac 32,\, |x_3|<1\}$ by
\begin{equation*}
\phi_{(1)}({\rm x}):=\left\{\begin{split}
%\sum_{i=1}^4c_i\phi\left(x_1,x_2,\frac{1-x_3}{i}\right)\quad&\mbox{for}\quad {\bf x}\in\overline{\Omega_{\rm ex}}\backslash\overline{\Omega_L}\cap\{x_3>1\},\\
\phi({\rm x})\quad&\mbox{for}\quad {\rm x}\in\overline{\Omega_L^{(1)}\cap \{|x_2|<1\}}(=\ol{\Om_L}),\\
\sum_{k=1}^5c_k\phi\left(x_1,{\rm sgn}\, x_2 (1+\frac 1k)-\frac{x_2}{k},x_3\right)\quad&\mbox{for}\quad {\rm x}\in\overline{\Omega_L^{(1)}\cap \{|x_2|>1\}},
\end{split}\right.
\end{equation*}
for constants $c_1,\cdots, c_5$ satisfying the equations
\begin{equation}
\label{interpolation}
\sum_{k=1}^5 c_k\left(-\frac{1}{k}\right)^m=1,\quad m=0,1,\ldots,4.
\end{equation}
It can be directly checked that  the linear system \eqref{interpolation} has a unique solution of $(c_1,c_2,c_3,c_4,c_5)=(15,-640,3645,-6144,3125)$.
Next, we extend $\phi_{(1)}$ onto $\Om_L^{\rm ex}:=\{{\rm x}\in \R^3: 0<x_1<L,\, |x_2|< \frac 32, |x_3|< \frac 32\}$ as follows:
\begin{equation*}
\phi_{{\rm ex}}({\rm x}):=\left\{\begin{split}
%\sum_{i=1}^4c_i\phi\left(x_1,x_2,\frac{1-x_3}{i}\right)\quad&\mbox{for}\quad {\bf x}\in\overline{\Omega_{\rm ex}}\backslash\overline{\Omega_L}\cap\{x_3>1\},\\
\phi_{(1)}({\rm x})\quad&\mbox{for}\quad {\rm x}\in\overline{\Om_L^{\rm ex}\cap \{|x_3|<1\}}(=\ol{\Om_L^{(1)}}),\\
\sum_{k=1}^5c_k\phi_{(1)}\left(x_1,x_2, {\rm sgn}\, x_3 (1+\frac 1k)-\frac{x_3}{k}\right)\quad&\mbox{for}\quad {\rm x}\in\overline{\Omega_L^{\rm ex}\cap \{|x_3|>1\}}.
\end{split}\right.
\end{equation*}
For each $r\in(0, \frac 14)$, let $\chi_r$ be a radially symmetric standard mollifier with a support in a ball of radius $r>0$. Then, a function $\phi^{(r)}$ given by
\begin{equation*}
\phi^{(r)}:=  \phi_{\rm ex}\ast \chi_r\quad\tx{in $\Om_L$}
\end{equation*}
provides a smooth approximation of $\phi$ in $\Om_L$, and we have $\displaystyle{\lim_{r\to 0+}}\|\phi-\phi^{(r)}\|_{H^4(\Om_L)}=0$. Furthermore, each $\phi^{(r)}$ satisfies the slip boundary condition
\begin{equation}\label{slip-bc-smooth-approx}
\der_{\bf n}\phi^{(r)}=0\quad\mbox{on}\quad\Gamma_w.
\end{equation}

For each $(\phi, W)\in \mcl{J}_{\eps, L}$, let $(\phi^{(r)}, W^{(r)})$ represent the smooth approximation of $(\phi, W)$ given by the procedure described in the above. And, for each $r\in(0, \frac 14)$, let $a_{ij}^{(r)}$ be given by
\begin{equation}
\label{definition-aij-r}
  a_{ij}^{(r)}:=a_{ij}(W^{(r)}, \nabla\phi^{(r)})\quad\tx{in $\Om_L$}.
\end{equation}
Then, from Lemmas \ref{lem-nei} and \ref{lemma-L1-hyperbolic-type}, the following lemma is directly given:
\begin{lemma}\label{lemma-smooth-coeff}
For the constant $\eps_2$ from Lemma \ref{lemma-L1-hyperbolic-type}, suppose that $\eps\le \frac{\eps_2}{2}$ in the definition of $\mcl{J}_{\eps, L}$ given in \eqref{Ite-set}. Then, there exists a small constant $\bar{r}>0$ depending on $(\gamma, J_0,  b_0, u_0, E_0, \delta)$ so that, for each $r\in(0, \bar{r}]$, the approximated coefficients $a_{ij}^{(r)}$ for $1\le i,j\le 3$ satisfy the following properties:
\begin{itemize}
\item[(a)]For the constant $\mu_1>0$ from Lemma \ref{lem-nei}, it holds that
\begin{equation*}
 \frac{\mu_1}{2}\le a_{11}^{(r)}, -a_{jj}^{(r)}\le \frac{2}{\mu_1}\quad\tx{in $\Om_L$}
 \end{equation*}
 for $j=2,3$. In fact, $a_{11}^{(r)}=1$;
\item[(b)] For $1\le i,j\le 3$, we have
$a_{ij}^{(r)}=a_{ji}^{(r)};$

\item[(c)]If $r\in(0, \bar{r}]$ is sufficiently small (depending on $\eps$), then it holds that
\begin{equation*}
      \sup_{1\le i,j\le 3}\|a_{ij}^{(r)}-\bar{a}_{ij}\|_{H^3(\Om_L)}
      \le C\eps
    \end{equation*}
    for a constant $C>0$ depending only on $(\gamma, J_0,  b_0, u_0, E_0, \delta)$. Moreover, the matrix $-[a_{ij}^{(r)}]_{i,j=2}^3$ is positive definite in $\Om_L$ with
    \begin{equation*}
      \frac{\mu_1}{4}\le -[a_{ij}^{(r)}]_{i,j=2}^3\le \frac{4}{\mu_1}\quad\tx{in $\Om_L$}
    \end{equation*}
    for the constant $\mu_1>0$ from Lemma \ref{lem-nei}.
\end{itemize}
\end{lemma}

By adjusting the smooth extension argument in the above, we can also approximate $(f_1, f_2)$ by smooth functions $(f_1^{(r)}, f_2^{(r)})$ so that we have
\begin{equation}
\label{compatibility-approx}
  \der_{\bf n} f_k^{(r)}=0\quad\tx{on $\Gam_w$},
\end{equation}
and $\displaystyle{\lim_{r\to 0+}\|f_k-f_k^{(r)}\|_{H^3(\Om_L)}=0}$
for $k=1,2$.

\subsubsection{Approximated linear boundary value problem and its a priori $H^1$ estimate} \label{subsubsection-H1}
Suppose that $\eps\le \frac{\eps_2}{2}$ in \eqref{Ite-set} so that Lemma \ref{lemma-smooth-coeff} holds, and let us fix $r\in(0, \bar{r}]$. For $a_{ij}^{(r)}$ given by \eqref{definition-aij-r}, we define a linear operator $\mcl{L}_1^{(r)}$ by
\begin{equation*}
\mcl{L}_1^{(r)}(\psi, \Psi):=  \sum_{i,j=1}^3a_{ij}^{(r)}\partial_{ij}\psi
+\bar{a}_1\partial_1\psi+\bar{b}_1\partial_1\Psi+\bar{b}_2 \Psi.
\end{equation*}
Assume that $\psi$ and $\Psi$ are smooth in $\ol{\Om_L}$, and that they solve the following linear boundary value problem:
\begin{equation}\label{lin-bd-simp-smooth}
\left\{\begin{split}
\mathcal{L}_1^{(r)}(\psi,\Psi)
=f_1^{(r)},\quad
\mathcal{L}_2(\psi,{\Psi})=f_2^{(r)}\quad&\mbox{in}\quad\Omega_L,\\
(\psi,\partial_1\psi)=(0,g),\quad\partial_1{\Psi}=0\quad&\mbox{on}\quad\Gamma_0,\\
{\Psi}=0\quad&\mbox{on}\quad\Gamma_L,\\
\der_{\bf n}\psi=0,\quad \der_{\bf n}{\Psi}=0\quad &\mbox{on}\quad\Gamma_w.
\end{split}\right.
\end{equation}
For the rest of Section \ref{subsection-more-details}, we let $a_{ij}$ and $f_k$ denote $a_{ij}^{(r)}$ and $f_k^{(r)}$, respectively for $i,j=1,2,3$, and $k=1,2$ unless otherwise specified.

Now we write down a detailed computation to get a priori $H^1$ estimate of $(\psi, \Psi)$. We emphasize that the focus of this paper is to show how to utilize the structure of the coefficients $a_{ij}^{(r)}$ and the slip boundary conditions for $\phi$ (an element from the iteration set $\mcl{J}_{\eps, L}$) and for $\psi$ (a part of solution to \eqref{lin-bd-simp-smooth}) so that the well-posedness of the boundary value problem \eqref{lin-bd-simp-smooth} can be achieved in a cylindrical domain of any dimension.
For a smooth function $\mcl{W}(x_1)$, which we call {\emph{an energy weight}}, let us define $I_{\rm H}(\psi, \Psi)$ by
\begin{equation}
\label{definition-Ih}
\begin{split}
  I_{\rm H}(\psi, \Psi)
  :=&\int_{\Om_L} \mcl{W}\psi_1\mcl{L}_1^{(r)}(\psi, \Psi)\,d\rx\\
  =&\int_{\Om_L} \sum_{i,j=1}^3 a_{ij} \psi_{ij} \mcl{W}\psi_1 \,d\rx+
  \int_{\Om_L} (\bar{a}_1\psi_1+\bar{b}_1\Psi_1+\bar{b}_2 \Psi) \mcl{W}\psi_1\,d\rx\\
  =:&I_{\rm H}^{p}+I_{\rm H}^{lot}.
  \end{split}
\end{equation}
Next, we take a closer look at the principal part $I_{\rm H}^{p}$ of $I_{\rm H}(\psi, \Psi)$. We again split $I_{\rm H}^{p}$ as follows:
\begin{equation}
\label{Ih}
  I_{\rm H}^{p}=\sum_{i=1}^3 \int_{\Om_L} \sum_{j=1}^3 a_{ij}\psi_{ij} \mcl{W}\psi_1 \,d\rx=:\sum_{i=1}^3 P_i.
\end{equation}
By the divergence theorem and the boundary condition $\psi_1=g$ on $\Gam_0$, we have
\begin{equation*}
\begin{split}
  P_1=&\int_{\Om_L} \sum_{j=1}^3 a_{1j}\mcl{W} \der_j\left(\frac{\psi_1^2}{2}\right) \,d\rx\\
  =&\int_{\Om_L}\nabla \cdot \left(a_{1j} \mcl{W} \frac{\psi_1^2}{2}\right)
  -\frac{\psi_1^2}{2}\nabla\cdot (\mcl{W} a_{1j})\,d\rx\\
  =&\int_{\Gam_L} a_{11}\mcl{W} \frac{\psi_1^2}{2}\,dS-\int_{\Gam_0}a_{11}\mcl{W} \frac{g^2}{2}\,dS
  -\int_{\Gam_w} \mcl{W} \frac{\psi_1^2}{2} (a_{1j})\cdot{\bf n}\,dS
  -\int_{\Om_L}\frac{\psi_1^2}{2}\nabla\cdot (\mcl{W} a_{1j})\,d\rx.
  \end{split}
\end{equation*}
Moreover, it directly follows from \eqref{def-aij}, \eqref{slip-bc-smooth-approx} and \eqref{definition-aij-r} that
\begin{equation*}
  (a_{1j})\cdot {\bf n}=\frac{(\der_1\vphi_0+\phi_1^{(r)} )(\nabla\vphi_0+ \nabla\phi^{(r)})\cdot {\bf n}}{\til{c}^2(\Phi_0+W^{(r)}, \nabla\vphi_0+\nabla\phi^{(r)})-(\der_1\vphi_0+\phi_1^{(r)})^2}=0\quad\tx{on $\Gam_w$},
\end{equation*}
and this yields that
\begin{equation}
\label{P1}
  P_1=\int_{\Gam_L} a_{11}\mcl{W} \frac{\psi_1^2}{2}\,dS-\int_{\Gam_0}a_{11}\mcl{W} \frac{g^2}{2}\,dS
  -\int_{\Om_L}\frac{\psi_1^2}{2}\nabla\cdot (\mcl{W} a_{1j})\,d\rx.
\end{equation}
To compute $P_2+P_3$, one needs more endeavor. Fix $i=2,3$, and let us consider $P_i$. By regarding $(a_{ij}\mcl{W}\psi_1)_{j=1}^3$ as a three dimensional vector field, we have
\begin{equation*}
  P_i=\int_{\Om_L} (a_{ij}\mcl{W}\psi_1)\cdot \nabla(\psi_i)\,d\rx.
\end{equation*}
Similar to the case of $P_1$, we integrate by parts and apply the divergence theorem to get
\begin{equation*}
\begin{split}
P_i=&\int_{\der \Om_L} (\sum_{j=1}^3 a_{ij}\hat{\bf e}_j)\cdot {\bm \nu}_{\rm out} \mcl{W}\psi_1 \psi_i\,dS
-\int_{\Om_L}\sum_{j=1}^3\der_j(a_{ij} \mcl{W}\psi_1)\psi_i \,d\rx\\
&=:P_i^{\rm bd}-R_i
\end{split}
\end{equation*}
for the outward unit normal ${\bm\nu}_{\rm out}$ on $\der \Om_L$.
By \eqref{def-aij}, \eqref{slip-bc-smooth-approx} and \eqref{definition-aij-r}, we have
\begin{equation}
\label{slip-bc-coeff}
  (\sum_{j=1}^3 a_{ij}\hat{\bf e}_j)\cdot {\bm\nu}_{\rm out}=  \frac{\til{c}^2(\Phi_0+W^{(r)}, \nabla\vphi_0+\nabla\phi^{(r)}) {\bf n}\cdot \hat{\bf e}_i}{\til{c}^2(\Phi_0+W^{(r)}, \nabla\vphi_0+\nabla\phi^{(r)})-(\der_1\vphi_0+\phi_1^{(r)})^2}\quad\tx{on $\Gam_w$}.
\end{equation}
Since ${\bf n}\cdot \hat{\bf e}_1=0$ on $\Gam_w$, we obtain that
\begin{equation*}
\sum_{i=2}^3  \frac{\til{c}^2 {\bf n}\cdot \hat{\bf e}_i}{\til{c}^2-(\der_1\vphi_0+\phi_1^{(r)})^2} \mcl{W}\psi_1\psi_i
=\frac{\til{c}^2 \mcl{W}\psi_1\nabla\psi\cdot {\bf n} }{\til{c}^2-(\der_1\vphi_0+\phi_1^{(r)})^2}
=0\quad\tx{on $\Gam_w$}
\end{equation*}
for $\til{c}=\til c(\Phi_0+W^{(r)}, \nabla\vphi_0+\nabla\phi^{(r)})$,
and this yields that
\begin{equation*}
  \sum_{i=2}^3 P_i^{\rm bd}=\sum_{i=2}^3 \int_{\Gam_L} a_{i1}\mcl{W}\psi_1\psi_i \,dS
\end{equation*}
because $\psi_1=0$ holds on $\Gam_0$.

Now, we compute $R_2+R_3$.
First, for each $i=2,3$, we check that
\begin{equation*}
\begin{split}
  \der_1(a_{i1} \mcl{W}\psi_1)\psi_i
  =&\der_1(a_{i1}\mcl{W}\psi_1\psi_i)
  -a_{i1}\mcl{W}\left(\frac{\psi_1^2}{2}\right)_i\\
  =&\der_1\left(a_{i1}\mcl{W}\psi_1\psi_i\right)
  -\der_i\left(a_{i1}\mcl{W}\frac{\psi_1^2}{2}\right)
  +
  \der_i(a_{i1}\mcl{W})\frac{\psi_1^2}{2}.
  \end{split}
\end{equation*}
Then we integrate by parts to obtain that
\begin{equation}
\label{R1}
 \sum_{i=2}^3 \int_{\Om_L}\der_1(a_{i1} \mcl{W}\psi_1)\psi_i\,d\rx
 =\sum_{i=2}^3\int_{\Gam_L} a_{i1}\mcl{W}\psi_1\psi_i\,dS+\sum_{i=2}^3\int_{\Om_L} \der_i(a_{i1}\mcl{W})\frac{\psi_1^2}{2}\,d\rx
\end{equation}
by using the boundary condition $\psi=0$ on $\Gam_0$ and \eqref{slip-bc-coeff}.
Next, we apply Lemma \ref{lemma-smooth-coeff}(b) to obtain that
\begin{equation*}
  \begin{split}
  &\sum_{i,j=2}^3\int_{\Om_L} \der_j(a_{ij}\mcl{W}\psi_1)\psi_i\,d\rx\\
  &= \sum_{i,j=2}^3\int_{\Om_L}(a_{ij}\psi_{j1}
  + \der_j a_{ij}\psi_1)\mcl{W}\psi_i\,d\rx\\
  &=\sum_{i,j=2}^3\int_{\Om_L}
  \frac 12 (a_{ij}\mcl{W}\psi_i\psi_j)_1-\frac 12 \der_1(a_{ij}\mcl{W})\psi_i\psi_j
   + \der_j a_{ij}\mcl{W}\psi_1\psi_i\,d\rx.
  \end{split}
\end{equation*}
Due to the boundary condition $\psi=0$ on $\Gam_0$, we get that
\begin{equation}
\label{R2}
\begin{split}
 &\sum_{i,j=2}^3\int_{\Om_L} \der_j(a_{ij}\mcl{W}\psi_1)\psi_i\,d\rx\\
 &=\sum_{i,j=2}^3\frac 12\int_{\Gam_L} a_{ij}\mcl{W}\psi_i\psi_j\,dS
 +\sum_{i,j=2}^3\int_{\Om_L}
   \der_j a_{ij}\mcl{W}\psi_1\psi_i
   -\frac 12 \der_1(a_{ij}\mcl{W})\psi_i\psi_j
   \,d\rx.
  \end{split}
\end{equation}
We combine \eqref{R1} and \eqref{R2} to have
\begin{equation*}
\begin{split}
  &R_2+R_3\\
  &=\sum_{i=2}^3\int_{\Gam_L} a_{i1}\mcl{W}\psi_1\psi_i\,dS+\sum_{i=2}^3\int_{\Om_L} \der_i(a_{i1}\mcl{W})\frac{\psi_1^2}{2}\,d\rx\\
  &\phantom{=}+\sum_{i,j=2}^3\frac 12\int_{\Gam_L} a_{ij}\mcl{W}\psi_i\psi_j\,dS
 +\sum_{i,j=2}^3\int_{\Om_L}
   \der_j a_{ij}\mcl{W}\psi_1\psi_i
   -\frac 12 \der_1(a_{ij}\mcl{W})\psi_i\psi_j
   \,d\rx,
  \end{split}
\end{equation*}
so we finally obtain that
\begin{equation}
\label{P2}
\begin{split}
  &P_2+P_3\\
  &=-\sum_{i=2}^3\int_{\Om_L} \der_i(a_{i1}\mcl{W})\frac{\psi_1^2}{2}\,d\rx\\
  &\phantom{=}-\sum_{i,j=2}^3\frac 12\int_{\Gam_L} a_{ij}\mcl{W}\psi_i\psi_j\,dS
 -\sum_{i,j=2}^3\int_{\Om_L}
   \der_j a_{ij}\mcl{W}\psi_1\psi_i
   -\frac 12 \der_1(a_{ij}\mcl{W})\psi_i\psi_j
   \,d\rx.
  \end{split}
\end{equation}
Note that $a_{11}\equiv 1$ in $\Om_L$. By \eqref{P1} and \eqref{P2}, the term $I_{\rm H}$, given by \eqref{definition-Ih}, is written as
\begin{equation}
\label{Ih-final}
\begin{split}
  I_{\rm H}(\psi, \Psi)
  =&\int_{\Gam_L} ( \psi_1^2-\sum_{i,j=2}^3 a_{ij}\psi_i\psi_j) \frac{\mcl{W}}{2}\,dS-\int_{\Gam_0}\mcl{W} \frac{g^2}{2}\,dS\\
  & +\int_{\Om_L}\left( \bar a_1 \mcl{W}-\frac 12\der_1\mcl{W} \right)\psi_1^2
  \,d\rx
  +\sum_{i,j=2}^3\int_{\Om_L}\frac 12 \der_1(a_{ij}\mcl{W})\psi_i\psi_j\,d\rx\\
  &+ \int_{\Om_L} (\bar{b}_1\Psi_1+\bar{b}_2 \Psi) \mcl{W}\psi_1\,d\rx-\sum_{i=2}^3\int_{\Om_L} \der_ia_{i1}\mcl{W}\psi_1^2\,d\rx
 -\sum_{i,j=2}^3\int_{\Om_L}
   \der_j a_{ij}\mcl{W}\psi_1\psi_i
   \,d\rx.
  \end{split}
\end{equation}
One can directly observe that the computations given in the above to obtain \eqref{Ih-final} can be easily adjusted for a cylindrical domain of any dimension with an arbitrary cross section. In fact, we do not require any condition for the cross section of a domain as long as its boundary is Lipscthiz continuous so that the divergence theorem can be applied.
Now that \eqref{Ih-final} is given, one can repeat the argument in \cite[Proof of Proposition 2.4]{BDXX} to establish a priori $H^1$ estimate of $(\psi, \Psi)$. Note that this is the essential part to prove Proposition \ref{proposition-main}.
Let us define $J_{\rm E}(\psi, \Psi)$ by
\begin{equation*}
\begin{split}
  J_{\rm E}(\psi, \Psi):=&-\int_{\Om_L} \Psi \mcl{L}_2(\psi, \Psi)\,d\rx\\
  =&\int_{\Om_L} |\nabla\Psi|^2+\bar{h}_1\Psi^2+\bar{h}_2\psi_1\Psi\,d\rx.
  \end{split}
\end{equation*}
Then we have
\begin{equation*}
  I_{\rm H}(\psi, \Psi)+J_{\rm E}(\psi, \Psi)=
  \int_{\Om_L} f_1 \mcl{W}\psi_1-f_2\Psi\,d\rx
\end{equation*}
which is equivalent to the following:
\begin{equation}
\label{energy-final}
  \begin{split}
  &\int_{\Gam_L} ( \psi_1^2-\sum_{i,j=2}^3 a_{ij}\psi_i\psi_j) \frac{\mcl{W}}{2}\,dS\\
  & +\int_{\Om_L}\left( \bar a_1 \mcl{W}-\frac 12\der_1\mcl{W} \right)\psi_1^2+\sum_{i,j=2}^3\int_{\Om_L}\frac 12 \der_1(a_{ij}\mcl{W})\psi_i\psi_j
  +|\nabla\Psi|^2+\bar{h}_1\Psi^2
  \,d\rx\\
  &+\int_{\Om_L} (\bar{b}_1\Psi_1+\bar{b}_2 \Psi) \mcl{W}\psi_1+\bar{h}_2 \psi_1\Psi\,d\rx-\sum_{i=2}^3\int_{\Om_L} \der_ia_{i1}\mcl{W}\psi_1^2\,d\rx
 -\sum_{i,j=2}^3\int_{\Om_L}
   \der_j a_{ij}\mcl{W}\psi_1\psi_i
   \,d\rx\\
   =&\int_{\Gam_0}\mcl{W} \frac{g^2}{2}\,dS+\int_{\Om_L} f_1 \mcl{W}\psi_1-f_2\Psi\,d\rx.
  \end{split}
\end{equation}

Now we explain how to find {\emph{the energy weight}} $\mcl{W}$ so that we have
\begin{equation}
\label{coercivity}
  \mbox{The left-hand side of \eqref{energy-final}}\ge \lambda_0 \int_{\Om_L} |\nabla\psi|^2+|\nabla\Psi|^2+\Psi^2\,d\rx
\end{equation}
for some constant $\lambda_0>0$.

Note that each $a_{ij}$ is close to $\bar{a}_{ij}$ provided that the constant $\eps$ in the definition of $\mcl{J}_{\eps, L}$ given in \eqref{Ite-set} and the constant $\bar r$ in Lemma \ref{lemma-smooth-coeff} are fixed sufficiently small. Let us first assume that
\begin{equation}
\label{assumption-coeff}
  a_{ij}=\bar{a}_{ij} \quad\tx{in $\Om_L$}.
\end{equation}
And, we apply Cauchy-Schwarz inequality to estimate {\emph{the coupling term }} $(\bar{b}_1\Psi_1+\bar{b}_2 \Psi) \mcl{W}\psi_1+\bar{h}_2 \Psi\psi_1$ as
\begin{equation*}
  |(\bar{b}_1\Psi_1+\bar{b}_2 \Psi) \mcl{W}\psi_1+\bar{h}_2 \Psi\psi_1|
  \le \frac 18 \Psi_1^2+\frac{\kappa_0}{4}\Psi^2+q_0(\mcl{W})\frac{\psi_1^2}{2}
\end{equation*}
with
\begin{equation*}%\label{eq-q4}
\begin{split}
\kappa_0:=\inf_{x_1\in[0,L_1]}\bar{h}_1(x_1),\quad
q_0(\mathcal{W})
:=4\left(\left(\bar{b}_1^2
+\frac{\bar{b}_2^2}{\kappa_0}\right)\mathcal{W}^2
+\frac{\bar{h}_2^2}{\kappa_0}\right)
\end{split}
\end{equation*}
for the constant $L_1$ from Lemma \ref{one-sol}. Note that it can be directly checked from \eqref{def-aij} and Lemma \ref{one-sol} that $\kappa_0>0$.
So we get
\begin{equation*}
\begin{split}
  &\mbox{The left-hand side of \eqref{energy-final} under the assumption of \eqref{assumption-coeff}}\\
  &\ge
  \int_{\Gam_L} ( \psi_1^2-\sum_{i,j=2}^3 \bar a_{ij}\psi_i\psi_j) \frac{\mcl{W}}{2}\,dS\\
  &\phantom{\ge} +\int_{\Om_L}\left( \bar a_1 \mcl{W}-\frac 12\mcl{W}' -\frac 12 q_0(\mcl{W}) \right)\psi_1^2+ \sum_{i=2}^3 \frac 12 \der_1(\bar a_{22}\mcl{W})\psi_i^2\,d\rx\\
  &\phantom{\ge} +\frac 34\int_{\Om_L}
|\nabla\Psi|^2+\bar{h}_1\Psi^2
  \,d\rx
  \end{split}
\end{equation*}
where we use the fact that $\bar{a}_{ij}=0$ for $i\neq j$ and $\bar{a}_{22}=\bar{a}_{33}$ in $\Om_L$.
Then we require the following conditions for $\mcl{W}$ to satisfy:
\begin{itemize}
\item[(i)] $\displaystyle{\min_{x_1\in[0, L]}}\mcl{W}(x_1)>0$;
\item[(ii)] $\displaystyle{\bar a_1 \mcl{W}-\frac 12\mcl{W}' -\frac 12 q_0(\mcl{W})>0}$;
\item[(iii)] $\displaystyle{ \frac 12 \der_1(\bar a_{22}\mcl{W})>0}$.
\end{itemize}
The existence of $\mcl{W}$ can be proved by setting up a differential equation for $\mcl{W}$ so that its solution satisfies all the three conditions stated in the above. Further details can be found in \cite[Step 3 in the proof of Proposition 2.4]{BDXX} so we state the following lemma without a proof.
\begin{lemma}
\label{lemma-existence-W}
Let us set $\mcl{A}(x_1)$ as
\begin{equation*}
 \mcl{A}(x_1):= \min\{\mcl{W}(x_1),
    \bar a_1(x_1) \mcl{W}(x_1)-\frac 12\mcl{W}'(x_1) -\frac 12 q_0(\mcl{W}(x_1)),
     \frac 12 \der_1(\bar a_{22}\mcl{W})(x_1)\}.
\end{equation*}
Given a small constant $\delta>0$, let the constant $L_1>0$ be from Lemma \ref{one-sol}. Then there exists a constant $L_1^*\in(0, L_1]$ depending on $\gam$, $J_0$, $u_0$, $E_0$, $b_0$, and $\delta$ so that one can find a smooth function $\mcl{W}$ on $[0, L_1^*]$ with satisfying the following properties:
\begin{itemize}
\item[(a)]$\displaystyle{
\min_{x_1\in[0, L_1^*]} \mcl{A}(x_1)\ge 0
}$;
\item[(b)] $\displaystyle{
\min_{x_1\in[0, L]} \mcl{A}(x_1)>0
}$ for any $L\in(0, L_1^*)$.
\end{itemize}
\end{lemma}

For the rest of the paper, we fix $L\in(0, L_1^*)$. Then Lemma \ref{lemma-existence-W} implies that there exists a constant $\lambda_0>0$ depending on $\gam$, $J_0$, $u_0$, $E_0$, $b_0$, $\delta$, and $L$ so that we have
\begin{equation*}
  \mbox{The left-hand side of \eqref{energy-final} under the assumption of \eqref{assumption-coeff}}
  \ge \lambda_0 \int_{\Om_L} |\nabla\psi|^2+|\nabla\Psi|^2+\Psi^2\,d\rx.
\end{equation*}
By Lemma \ref{lemma-smooth-coeff}, one can fix a small constant $\om_0>0$ depending on $\gam$, $J_0$, $u_0$, $E_0$, $b_0$, $\delta$, and $L$ so that if the constant $\eps$ in the definition of the iteration set $\mcl{J}_{\eps,L}$, given by \eqref{Ite-set}, satisfies
\begin{equation*}
     \eps\le \om_0,
   \end{equation*}
then any smooth solution $(\psi, \Psi)$ to the linear boundary value problem \eqref{lin-bd-simp-smooth} associated with $(\phi, W)\in \mcl{J}_{\eps,L}$ satisfies the estimate
\begin{equation*}
  \mbox{The left-hand side of \eqref{energy-final}} \ge \frac{\lambda_0}{2} \int_{\Om_L} |\nabla\psi|^2+|\nabla\Psi|^2+\Psi^2\,d\rx,
\end{equation*}
from which the following a priori $H^1$ estimate of $(\psi, \Psi)$ is obtained:
\begin{equation}
\label{apriori-H1}
  \|\psi\|_{H^1(\Om_L)}+\|\Psi\|_{H^1(\Om_L)}
  \le C\left(\|g\|_{L^2(\Gam_0)}+\sum_{k=1}^2 \|f_k\|_{L^2(\Om_L)}\right)
\end{equation}
for a constant $C>0$ chosen depending on $\gam$, $J_0$, $u_0$, $E_0$, $b_0$, $\delta$, and $L$.

\subsubsection{A brief discussion on higher order derivative estimates of $(\psi, \Psi)$}\label{subsubsection-H4}
 Once a priori $H^1$ estimate of $(\psi, \Psi)$ is achieved, higher order derivative estimates are given by a bootstrap argument. First of all, we rewrite $\mathcal{L}_1^{(r)}(\psi,\Psi)
=f_1$ and $\mathcal{L}_2(\psi,{\Psi})=f_2$ as
\begin{equation}
\label{definition-operator-separated}
\begin{split}
&\mathfrak{L}_1\psi:=\sum_{i,j=1}^3 a_{ij}\psi_{ij}+\bar a_1\psi_1=f_1-\bar b_1\Psi_1-\bar b_2\Psi=:F_1,\\
&\mathfrak{L}_2\Psi:=\Delta \Psi-\bar h_1 \Psi=f_2+\bar h_2\psi_1=:F_2.
\end{split}
\end{equation}
We regard $\Psi$ as a solution to the linear boundary value problem
\begin{equation}
\label{lbvp-Psi}
\begin{split}
&\mathfrak{L}_2\Psi=F_2\quad \mbox{in $\Om_L$},\\
&\Psi_1=0\quad\mbox{on $\Gam_0$},\quad
\der_{\bf n}\Psi=0\quad\tx{on $\Gam_w$},\quad \Psi=0\quad\tx{on $\Gam_L$},
\end{split}
\end{equation}
and apply the method of reflection and the standard elliptic theory(a local $H^2$ estimate of weak solution to elliptic boundary value problem) to obtain that
\begin{equation}
\label{H2-Psi}
  \|\Psi\|_{H^2(\Om_L)}\le C(\|F_2\|_{L^2(\Om_L)}+\|\Psi\|_{L^2(\Om_L)}).
\end{equation}
Once the estimate \eqref{H2-Psi} is obtained, we move on to the boundary value problem
\begin{equation}
\label{lbvp-psi}
  \begin{split}
  &\mathfrak{L}_1\psi=F_1\quad\tx{in $\Om_L$},\\
  &(\psi, \psi_1)=(0, g)\quad\tx{on $\Gam_0$}, \quad \der_{\bf n}\psi=0\quad\tx{on $\Gam_w$}.
  \end{split}
\end{equation}
To estimate $H^2$-norm of $\psi$, we apply the method of  standard hyperbolic estimate. For each $t\in(0, L]$, let us define
\begin{equation*}
  \Om_t:=\Om_L\cap\{x_1<t\},\quad \Gam_t:=\ol{\Om_L}\cap\{x_1=t\},
\end{equation*}
and consider the following expression:
\begin{equation}
\label{H2-comp-psi}
  \int_{\Om_t} \psi_{11}\der_1(\mathfrak{L}_1\psi)\,d\rx
  =\int_{\Om_t} \der_1F_1 \psi_{11}\,d\rx.
\end{equation}
Let us set
\begin{equation*}
  \mathfrak{R}:=\|f_1\|_{H^1(\Om_L)}+\|f_2\|_{L^2(\Om_L)}
  +\|g\|_{C^1(\Gam_0)}.
\end{equation*}
Then, by  a lengthy but straightforward computation with a careful integration by parts and using \eqref{apriori-H1}, one can derive from \eqref{H2-comp-psi} that
\begin{equation}
\label{gronwall-pre}
  \int_{\Gam_t} |\nabla\psi_1|^2\,d\rx
  \le \alp \int_{\Om_t}|\nabla \psi_1|^2+|D^2_{(x_2,x_3)}\psi|^2\,d\rx
  +\beta \mathfrak{R}^2
\end{equation}
for constants $\alp>0$ and $\beta>0$ fixed depending only on $\gam$, $J_0$, $u_0$, $E_0$, $b_0$, $\delta$, and $L$. In order to estimate $\|D^2_{(x_2, x_3)}\psi\|_{L^2(\Om_t)}$, an additional step is necessary. %For any fixed $s\in(0,L)$, set
%\begin{equation*}
%  \Gam_s:=\{s\}\times \mcl{D} (\subset \Om_L).
%\end{equation*}
On $\Gam_s(\subset \R^2)$ for $s\in(0,L)$, $\psi$ satisfies
\begin{equation}
\label{lbvp-2d-eqn}
\begin{split}
  &\underset{(=\bar{a}_{22}\sum_{i=2}^3\psi_{ii})}{\underbrace{\sum_{i,j=2}^3 \bar a_{ij}\psi_{ij}}}=F_1-\bar a_1\psi_1-\sum_{{1\le i,j\le 3}\atop{i\,\tx{or}\,j=1}} a_{ij}\psi_{ij}+\sum_{i,j=2}^3 (\bar a_{ij}-a_{ij})\psi_{ij}=:\mathfrak{F}_1.
  %&\der_{\bf n}\psi=0\quad\tx{on $\{s\}\times \der \mcl{D}$}.
  \end{split}
\end{equation}
And, on $\der \Gam_s(=\{s\}\times \der \mcl{D})$, $\psi$ satisfies the slip boundary condition
\begin{equation}
\label{lbvp-2d-bc}
  \der_{\bf n}\psi=0.
\end{equation}
We regard \eqref{lbvp-2d-eqn} and \eqref{lbvp-2d-bc} as a two dimensional elliptic boundary value problem for $\psi$ on a square $\Gam_s$, and apply the method of reflection and the standard elliptic theory to obtain that
\begin{equation*}
  \|D^2_{(x_2, x_3)}\psi\|_{L^2(\Gam_s)}
  \le C_*(\|F_1\|_{L^2(\Gam_s)}+\|\psi_1\|_{L^2(\Gam_s)}
  +\|\nabla\psi_1\|_{L^2(\Gam_s)}+\eps \|D^2_{(x_2, x_3)}\psi\|_{L^2(\Gam_s)})
\end{equation*}
for a constant $C_*>0$ depending only on $\gam$, $J_0$, $u_0$, $E_0$, $b_0$, and $\delta$. So if $\eps>0$ is fixed small depending on $\gam$, $J_0$, $u_0$, $E_0$, $b_0$, and $\delta$, it follows from the estimate right above that
\begin{equation}
\label{rem-deriv-H2}
  \|D^2_{(x_2, x_3)}\psi\|^2_{L^2(\Gam_s)}
  \le 2C_*(\|F_1\|^2_{L^2(\Gam_s)}+\|\psi_1\|^2_{L^2(\Gam_s)}
  +\|\nabla\psi_1\|^2_{L^2(\Gam_s)}).
\end{equation}
By substituting this estimate into \eqref{gronwall-pre}, it is obtained that
\begin{equation*}
%\label{gronwall-main}
  \int_{\Gam_t} |\nabla\psi_1|^2\,d\rx
  \le \alp' \int_{\Om_t}|\nabla \psi_1|^2\,d\rx
  +\beta' \mathfrak{R}^2
\end{equation*}
for constants $\alp'>0$ and $\beta'>0$ fixed depending only on $\gam$, $J_0$, $u_0$, $E_0$, $b_0$, $\delta$, and $L$. And, this combined with Gronwall's inequality yields that
\begin{equation}
\label{H2-pre}
  \|\nabla\psi_1\|_{L^2(\Om_L)}\le C\mathfrak{R}
\end{equation}
for a constant $C>0$ depending only on $\gam$, $J_0$, $u_0$, $E_0$, $b_0$, $\delta$, and $L$. Next, we integrate the inequality \eqref{rem-deriv-H2} in $x_1$-direction over the interval $(0,L)$ and apply the estimate \eqref{H2-pre} to conclude that
\begin{equation}
\label{H2}
  \|\psi\|_{H^2(\Om_L)}\le C\mathfrak{R}.
\end{equation}

A priori $H^3$ or $H^4$ estimates can be obtained similarly. But there is one issue to be checked carefully for higher order derivative estimates. Note that $\Om_L$ is a Lipschitz domain because its cross section is a square. As we have already mentioned in the a priori $H^2$ estimate of $(\psi, \Psi)$ (even though no detail is given), the method of reflection is applied so that the corner points in the set $\{{\rm x}=(x_1, x_2, x_3)\in \der \Om_L: |x_2|=|x_3|=1 \}$ can be treated as interior points. Since the slip boundary conditions are prescribed for both $\psi$ and $\Psi$ on $\Gam_w$, we use even extensions of $\psi$ and $\Psi$ about either $|x_2|=1$ or $|x_3|=1$.

Thankfully, a priori $H^3$ estimate of $(\psi, \Psi)$ can be directly given by the method of reflection just as in a priori $H^2$ estimate. In order to get a priori $H^3$ estimate of $\Psi$, we need for the even extension of $F_2$  about $|x_2|=1$(or $|x_3|=1$) to be $H^1$ across $|x_2|=1$ (or $|x_3|=1$). In order to get a priori $H^3$ estimate of $\psi$, a similar condition is required but it is a bit more complicated. For each $s\in(0, L)$, the function $\mathfrak{F}_1(s,\cdot)$ given in \eqref{lbvp-2d-eqn} is defined on a square $\mcl{D}$, the cross section of $\Om_L$. As a function of $(x_2, x_3)$, if the even extension of $\mathfrak{F}_1(s,\cdot)$ about $|x_2|=1$ (or $|x_3|=1$) is $H^1$ across $|x_2|=1$(or $|x_3|=1$), then one can establish a priori estimate of $\|D^3_{(x_2, x_3)}\psi\|_{L^2(\Gam_s)}$. Once this estimate is achieved, a priori $H^3$ estimate of $\psi$ can be obtained by adjusting the computations given to get a priori $H^2$ estimate of $\psi$ although a few more extra steps may be needed.
Finally, it is not hard to check that an even extension of a $H^1$ function is still $H^1$ across the reflection boundary. So a priori $H^3$ estimate of $(\psi, \Psi)$ can be given without any complication.

A priori $H^4$ estimate of $(\psi, \Psi)$ requires a more endeavor though. It is same as $H^3$ estimate in that we apply the method of reflection to get a priori $H^4$ estimate of $\Psi$ by using the boundary value problem \eqref{lbvp-Psi}. For $H^4$ estimate, however, we additionally need for the even extension of $F_2$ about $|x_2|=1$(or $|x_3|=1$) to be $H^2$ across $|x_2|=1$ (or $|x_3|=1$). By using the compatibility condition of $f_2$ stated in \eqref{compatibility-approx}, the definition of $F_2$ given in \eqref{definition-operator-separated} and the slip boundary condition of $\psi$ on $\Gam_w$, one can directly check that $\der_{\bf n} F_2=0$ on $\Gam_w$. From this, it can be directly checked that the even extension of $F_2$ about $|x_2|=1$(or $|x_3|=1$) is $H^2$ across $|x_2|=1$(or $|x_3|=1$) thus a priori $H^4$ estimate of $\Psi$ can be achieved in $\Om_L$. Now, the remaining question is whether the even extension of $\mathfrak{F}_1(s, \cdot)$ about $|x_2|=1$(or $|x_3|=1$) is $H^2$ across $|x_2|=1$(or $|x_3|=1$) as a two variable function of $(x_2, x_3)$ for each $s\in(0, L)$. Alternatively, can we show that $\der_{\bf n}\mathfrak{F}_1=0$ on $\Gam_w$? If the answer is yes, then we can establish a priori estimate of $\|D^4_{(x_2, x_3)}\psi\|_{L^2(\Gam_s)}$ for each $s\in (0,L)$ thus establish a priori $H^4$ estimate of $\psi$ in $\Om_L$. Fortunately, the answer is yes. But, in order to verify it, a lengthy computation must be done. In this paper, we do not put any further details as the computation is very straightforward. Instead, we list two important points of the computation:
\begin{itemize}
\item[(i)] The compatibility condition \eqref{slip-bc-smooth-approx} and the particular structure of $a_{ij}^{(r)}$ given by \eqref{definition-aij-r} must be collaborated in a subtle way;
\item[(ii)] The slip boundary condition $\der_{\bf n}\psi=0$ on $\Gam_w$ must be used appropriately.
\end{itemize}

Finally, we note that all the arguments stated in the above can be extended into a cylindrical domain with a rectangular cross section in $\R^n$ for any $n\ge 3$. Also, they can be applied to a cylindrical domain if a cross section has a geometrically symmetric property such as an axisymmetric property for instance. This case will be studied in the forthcoming work. Of course, more important question is what happens if a cross section is in an arbitrary shape. Can we still prove the well-posedness of the boundary value problem \eqref{lin-bd} in $H^4$? At the moment, the answer is unknown. But this is certainly an important problem to be studied in the future.

\subsubsection{A brief description on how to establish the existence of a $H^4$ solution to the boundary value problem \eqref{lin-bd} } In \S \ref{subsubsection-H1}--\S\ref{subsubsection-H4}, we have explained how to establish a priori $H^4$-estimate of $(\psi, \Psi)$ by assuming that $(\psi, \Psi)$ is smooth in $\Om_L$ up to the boundary, and that it solves \eqref{lin-bd-simp-smooth}.
For readers who are familiar with the results in \cite{BDXX}, they would immediately see how to complete a proof of Proposition \ref{proposition-main} by using the argument given in \S \ref{subsubsection-H1}--\S\ref{subsubsection-H4}. To make this paper self-contained, we briefly explain how to complete a proof of Proposition \ref{proposition-main}. Overall, we take the following steps:\\
\quad\\
{\emph{Step 1:}} Fix $r>0$ sufficiently small. Given $(\phi, W)\in \mcl{J}_{\eps, L}$ and $(f_1, f_2)$, we approximate them as $(\phi^{(r)}, W^{(r)})$ and $(f_1^{(r)}, f_2^{(r)})$ by following the argument given in \S\ref{subsubsection-approximation}. Then, we set up an approximated linear boundary value problem \eqref{lin-bd-simp-smooth}.\\
\quad\\
{\emph{Step 2:}} (Galerkin approximation) Let $\mathfrak{e}:=\{\omega_k\}_{k=0}^{\infty}\subset \R_+$
be the collection of all eigenvalues to the following two dimensional eigenvalue problem
\begin{equation}\label{basis}
-\Delta_{{(x_2, x_3)}}\eta=\omega\eta\quad\mbox{in}\quad \mathcal{D}, \quad\der_{\bf n}\eta=0\quad\mbox{on}\quad\partial \mathcal{D},
\end{equation}
and let $\mathfrak{E}:=\{\eta_k\}_{k=0}^{\infty}\subset C^{\infty}(\mcl{D})$ be the collection of all eigenfunctions satisfying that $\displaystyle{-\Delta_{{(x_2, x_3)}}\eta_k=\omega_k\eta_k}$ for each $k=0, 1, 2, \cdots$. One can choose the set $\mathfrak{E}$ so that it forms an orthonormal basis in $L^2(\mathcal{D})$ and an orthogonal basis in $H^1(\mathcal{D})$.

Define the standard inner product $\langle\cdot,\cdot\rangle$ in $L^2(\mathcal{D})$ by
\begin{equation*}
\langle \zeta_1,\zeta_2\rangle=\int_{\mathcal{D}}\zeta_1({\rm x'})\zeta_2({\rm x'})d{\rm x'}\quad\mbox{for ${\rm x}'=(x_2, x_3)\in \mcl{D}$.}
\end{equation*}
Fix $m\in\mathbb{N}$. For ${\rm x}=(x_1,{\rm x'})\in\Omega_L$, let us set $(\psi_m,\Psi_m)$ as
\begin{equation*}%\label{VW-form}
\begin{split}
&\psi_{m}(x_1,{\rm x'}):=\sum_{j=0}^{m} \theta_j(x_1)\eta_j({\rm x'}),\\
&\Psi_{m}(x_1,{\rm x'}):=\sum_{j=0}^{m}\Theta_j(x_1)\eta_j({\rm x'}).
\end{split}
\end{equation*}
For each fixed $m\in \mathbb{N}$, we find $(\theta_j(x_1), \Theta_j(x_1))$ for $j=0,\cdots, m$ so that it satisfies
\begin{equation}\label{VWm}
\langle\mathcal{L}_1^{(r)}(\psi_m,\Psi_m),\eta_k\rangle=\langle {f}_1^{(r)},\eta_k\rangle,\quad
\langle\mathcal{L}_2(\psi_m,\Psi_m),\eta_k\rangle=\langle {f}_2^{(r)},\eta_k\rangle
\end{equation}
for all $k=0,1,\dots,m$, and
\begin{equation}\label{VW-b}
\left\{\begin{split}
\psi_m=0,\quad \partial_1 \psi_m=\sum_{j=0}^m\langle g,\eta_j\rangle\eta_j,\quad \partial_1 \Psi_m=0\quad&\mbox{on}\quad\Gamma_0,\\
\Psi_m=0\quad&\mbox{on}\quad\Gamma_L.
\end{split}\right.
\end{equation}
From \eqref{VWm}, a system of ordinary differential equations for $\{(\theta_j, \Theta_j)\}_{j=0}^m$ is given on the interval $(0, L)$. And, \eqref{VW-b} yields boundary conditions for $\{(\theta_j, \Theta_j)\}_{j=0}^m$ at $x_1=0$ and $L$. Since all the coefficients of the differential equations, $(f_1^{(r)}, f_2^{(r)})$ and the eigenfunctions $\{\eta_j\}_{j=0}^{\infty}$ are smooth, if the boundary value problem for $\{(\theta_j, \Theta_j)\}_{j=0}^m$ has a solution then it is smooth on $[0, L]$ thus $(\psi_m, \Psi_m)$ is smooth in $\ol{\Om_L}$. Most importantly, a priori $H^1$ estimate of a smooth solution to \eqref{lin-bd-simp-smooth}, discussed in \S \ref{subsubsection-H1}, is used in establishing the well-posedness of the boundary value problem for $\{(\theta_j, \Theta_j)\}_{j=0}^m$. The unique existence of smooth function $(\psi_m, \Psi_m)$ satisfying \eqref{VWm} and \eqref{VW-b} can be proved by combining the results given in \S \ref{subsubsection-H1}--\S\ref{subsubsection-H4} and Fredholm alternative theorem applied to the system of ODEs for $\{(\theta_j, \Theta_j)\}_{j=0}^m$.\\
\quad\\
{\emph{Step 3:}} Finally, we take the limit $m\to \infty$ of $\{(\psi_m, \Psi_m)\}_{m\in \mathbb{N}}$ in order to obtain a $H^4$ solution to the approximated boundary value problem \eqref{lin-bd-simp-smooth}, then take another limit $r\to 0$ to obtain a $H^4$ solution to the boundary value problem \eqref{lin-bd}. The uniqueness of a solution is directly given from \eqref{hat-est}, thus Proposition \ref{proposition-main} is proved.

In the three steps described in the above, {\emph{Step 2}} is obviously the most important step. The approach using the Galerkin's approximation may seem as a standard method. But we must emphasize that this approach is possible thanks to the slip boundary conditions $\der_{\bf n}\psi=0$ and $\der_{\bf n}\Psi$=0 on $\Gam_w$, given in \eqref{lin-bd} because the same boundary conditions prescribed for $\psi$ and $\Psi$ enable us to describe $(\psi, \Psi)$ as the $H^4$ limit of their finite dimensional approximations in the eigenfunction space $\mathfrak{E}$. Also, we remark that the compatibility condition $\der_{\bf n}g=0\quad\tx{on $\ol{\Gam_0}\cap \ol{\Gam_w}$}$(that corresponds to the compatibility condition for $u_{\rm en}$ stated in \eqref{223}), which we have never mentioned so far, plays a crucial role in {\emph{Step 2}}.

Finally, we explain how the compatibility conditions \eqref{223} given for $(b, E_{\rm en}, \Phi_{\rm ex})$ are used. One may notice that the nonlinear boundary value problem \eqref{nlbvp-perturbation} contains the nonhomogeneous boundary conditions for $\Psi$, and so does its approximated linear boundary value problem \eqref{lin-pro2}. The boundary conditions become homogeneous by replacing $\Psi$ with $\Psi-(x_1-L)g_2(x_2, x_3)+\Psi_{\rm ex}$, therefore all of $(b, E_{\rm en}, \Phi_{\rm ex})$ are part of the nonhomogeneous terms of $(f_1,f_2)$. Then the compatibility conditions \eqref{223} given for $(b, E_{\rm en}, \Phi_{\rm ex})$ are contributed to satisfy the compatibility conditions \eqref{compatibility-original}.

%%%%%%%%%%%%%%%%%%%%%%%%%%%%%%%%%
%%%%%%%%%%%%%%%%%%%%%%%%%%%%%%%%%%%%

\vspace{.25in}
\noindent
{\bf Acknowledgements:}
The research of Myoungjean Bae was supported in part by  Samsung Science and Technology Foundation under Project Number SSTF-BA1502-51.
The research of Hyangdong Park was supported in part by the National Research Foundation of Korea (NRF) grant funded by the Korea government (MSIT and MOE) (No. 2015R1A5A1009350 and  No. 2020R1I1A1A01058480).

\bigskip
\bibliographystyle{siam}
%\bibliography{References}

\end{document}